# TRANSPOSE ON VERTEX SYMMETRIC DIGRAPHS

*Vance Faber*
*Center for Computing Sciences, Institute for Defense Analyses*

**Abstract**. *In [2] (and earlier in [3]), we defined several global communication tasks (universal exchange, universal broadcast, universal summation) on directed or vertex symmetric networks. In that paper, we mainly focused on universal broadcast. In [4], we discussed universal sum. Here we continue this work and focus on universal exchange (transpose) in directed vertex symmetric graphs.*

**Introduction**. In universal exchange, processor *i* has a vector of data $a_{ij}$, where the element $a_{ij}$ is a packet needed by processor *j*. Thus we start with a matrix $A$ where each processor has a separate row and we end the task with each processor having a separate column, in other words, we are computing the transpose.

This gives rise to two fundamental questions:

*Given a directed vertex symmetric graph $G$ (or any directed graph for that matter), how many time steps $\tau(G)$ does it take to perform a transpose?*

*Given two directed vertex symmetric graphs (or general directed graphs for that matter), how do we compare them with respect to their ability to perform a transpose?*

**Graph symmetries**. One important property of vertex symmetric graphs is that the number of edges directed into and out of a given vertex is a constant *d*. In our model of communication, we assume that on every time step all the edges in the graph can be used simultaneously, that is, each processor can simultaneously exchange a single packet of data with all of its neighbors.

From here on out, we consider only directed Cayley coset graphs (equivalently, directed vertex symmetric graphs) $G = (\Gamma, \Delta, H)$. Briefly, a Cayley coset graph $G$ has $P = [\Gamma : H]$ cosets of the subgroup $H$ in the group $\Gamma$ as vertices and degree $d = |\Delta|$ where $\Delta$ is a collection of elements in $\Gamma$ and $H \cup \Delta$ generates $\Gamma$. (Note that we allow multi-edges here but we insist that the graph be connected.) The edges are given by $(gH, g\delta H)$ for $\delta \in \Delta$. A necessary and sufficient condition for the edges to be well-defined between cosets in this way is that the union of the left and right cosets $\Delta H = H\Delta$. In the case where $H$ is the identity subgroup, $G$ is a directed Cayley graph. Let the number of vertices of distance $k$ from a given vertex be $n_k$. A simple fundamental lower bound on $\tau(G)$ is given by the following lemma.

*Lemma 1. The time for transpose is at least as great as the average diameter $\theta(G)$ of the graph:*

$$\tau(G) \geq \theta(G) = \left\lceil \sum k n_k / d \right\rceil. \qquad \text{(Equation 1)}$$

*Proof.* Each element $a_{ij}$ travels on a path (task graph) $T_{ij}$ from $i$ to $j$. Since the graph is vertex symmetric, the total length of these paths from a single vertex is at least $\sum k n_k$ and that means that all these transfers occupy at least $P \sum k n_k$ edges. If we could use all $Pd$ edges at once on each time step, we would achieve a time equal to $\theta(G)$.

**Problems achieving the optimum value**. In order to give an actual value for the transpose, we have to produce a communication graph, a schedule that lists for each time which edges transfer which data. This was discussed in [2, Definitions 2.2 and 2.13]. In the case of transpose, this is a collection of directed paths $T_{ij}$ in $G$ for each ordered pair of vertices with the properties that the times increase on each $T_{ij}$ and no edge is labeled twice with the same time. (Note that this can also be thought of as embedding the complete directed graph into $G$.) If we want to achieve $\theta(G)$, the paths $T_{ij}$ must be nearly the shortest paths from $i$ to $j$. (There is a bit of leeway in the case where $\theta(G)$ has been rounded up.) One problem is that any collection of shortest paths might have some edges represented more often than others. Another problem is that shortest paths from one vertex might conflict with shortest paths from another vertex. We discuss these problems one at a time.

**Paths without conflicts**. Suppose that $G = (\Gamma, \Delta, H)$ is a Cayley graph ($H$ is the identity subgroup.) Suppose we express each element $g$ in $G$ by a word $w(g)$ in (in other words, a product of) the members of $\Delta$. Then $w$ identifies a path from the identity $e$ to the element $g$ in $G$. Each entry in $w$ can be thought of as a single directed edge in $G$. Let $\delta(j)$ be the generator in entry $j$ of this path.

*Theorem 2. Suppose we label the edges of the paths $w(g)$ with increasing times so that no generator appears more than once with a given time. Then the collection of paths $hw(g)$ for all $h$ in $G$ is a communication graph for transpose on $G$ where the labels on the edges of the $hw(g)$ are identical to the labels on $w(g)$.*

*Proof.* We have to show that no edge has the same time twice. Suppose the same edge appears in $h_1 w(g_1)$ and $h_2 w(g_2)$ with the same time. Thus this edge can be represented in two ways:

$(h_1 g, h_1 g \delta(j))$ and $(h_2 h, h_2 h \delta(k))$.

So we must have

$h_1 g = h_2 h$

and

$h_1 g \delta(j) = h_2 h \delta(k)$

for some group elements $g$ and $h$. This implies $\delta(j) = \delta(k) = \delta$. Now the times on the edge $(h_1 g, h_1 g \delta)$ in the path $h_1 w(g_1)$ and the edge $(h_2 h, h_2 h \delta)$ in the path $h_2 w(g_2)$ are the same as the times on the edge $(g, g\delta)$ in the word $w(g_1)$ and the edge $(h, h\delta)$ in the word $w(g_2)$, respectively, by construction. Unless $g = h$ and $w(g_1) = w(g_2)$, this violates our assumption that no two edges with the same generator in the set of paths can have the same time.

This theorem allows us to schedule a transpose on $G$ by creating a consistent set of times on the paths from the identity to the other elements. Once we have chosen a set of such paths (words) $W$, we can calculate a lower bound on $\tau(G)$ with respect to these paths.

*Definition 1.* We let the *regular bound* on transpose be

$$\psi(G) = \min_W \max_\delta |\{\delta : \delta(j) \in W, \delta(j) = \delta\}|.$$

*Remark.* Since we have forced each occurrence of the generator $\delta$ to have a different time, the time that we obtain for transpose in this way is at least as great as $\psi(G)$. Also

$$\sum k n_k / d = \min_W \sum_\delta |\{\delta : \delta(j) \in W, \delta(j) = \delta\}| / d \leq \psi(G)$$

so $\theta(G) \leq \psi(G)$. What is not clear is whether there is any provable relationship between the regular bound and the optimal transpose. The ideal situation would be if each generator appears roughly the same number of times in the optimal $W$ but often that cannot be arranged without increasing the average size of the words in $W$ beyond the average diameter. There is no obvious reason why the regular bound is a lower bound on the time for a transpose on $G$ because there is no necessity for times to be consistently associated with generators throughout the graph so the same time could appear twice on the same generator in the paths from a fixed vertex to the other vertices. In addition, for $\psi(G)$ to be the maximum time in the labeling, every time must be assigned to one of the occurrences of the generator that appears most often. This might not be possible.

*Question 1.* Is $\tau(G) = \psi(G)$?

**Spanning factorizations**. The difficulty in routing on vertex symmetric graphs that are not Cayley graphs is the failure to be able to label edges with specific generators. As an example, the problem is seen in the Petersen graph [7] that is vertex symmetric but cannot be represented as a Cayley graph.

The Petersen graph is the complement of the line graph of $K_5$; the vertices can be represented as 2-element subsets and $\{u,v\}$ and $\{x,y\}$ are connected if they are disjoint. The automorphism group is isomorphic to $S_5$; each automorphism $\alpha$ is a derived permutation of vertices of $K_5$. (That is, given an automorphism $\alpha$ it is a map $\alpha\{u,v\} = \{\hat{\alpha}(u), \hat{\alpha}(v)\}$ for a permutation $\hat{\alpha}$.) The stabilizer subgroup $B$ of the vertex $\{1,2\}$ is thus the set of twelve derived permutations that fix the set $\{1,2\}$. One set of generators consists of the permutations $\delta_1$, $\delta_2$, $\delta_3$ derived from $(13)(24)$, $(13)(25)$, $(14)(25)$. Let $x = (34)(125)$. Then $xB \neq B$ are two distinct vertices of $G$. Also $(\delta_1)^{-1} x \delta_1 = (12)(345) \in B$ so $x\delta_1 B = \delta_1 B$. This means that there are two distinct edges going out from the vertex $x\delta_1 B = \delta_1 B$, one to $xB$ and the other to $B$, that can be thought to have the label $\delta_1$ depending on which coset representative is used, $x\delta_1$ or $\delta_1$.

We now discuss a possible solution to this problem. This section is taken from [1]. A 1-factor of a directed graph $G$ is a subgraph with both in-degree and out-degree equal to 1. (Some authors have called this a 2-factor. Our definition seems more consistent with the notation in undirected graphs. For example, if the edges are all bi-directional and the factor is a union of 2-cycles, then this would be an ordinary 1-factor in an undirected graph.)  It is known (Petersen 1891, see [7]) that every regular directed graph with in-degree and out-degree $d$ has a 1-factoring with $d$ 1-factors.  For completeness, we give the proof here.

*Fact 1.  Every directed graph $G$ where the in-degree and out-degree of every vertex is $d$ has a edge disjoint decomposition into $d$ 1-factors.*

*Proof.*  Form an auxiliary graph $B$ with two new vertices $u'$ and $u''$ for each vertex $u$.  The edges of $B$ are the pairs $(u', v'')$ where $(u, v)$ is a directed edge in $G$.  The undirected graph $B$ is bipartite and regular with degree $d$ so it can be decomposed into $d$ 1-factors.  Each of these 1-factors corresponds to a directed 1-factor in $G$.

In order to create global routing schemes on $G$, we consider regular graphs with factorizations with additional properties.

*Definition.* Let $F_1$, $F_2$, $\cdots$, $F_d$ be the factors in a 1-factoring of $G$. We call a finite string of symbols from the set $\{F_i\}$ a *word*. If $v$ is a vertex and $\omega$ is a word, then $v\omega$ denotes the directed path (and its endpoint) in $G$ starting at $v$ and proceeding along the unique edge corresponding to each consecutive factor represented in the word $\omega$.  If $G$ is a graph with $n$ vertices, we say that a 1-factoring and a set of $n$ words $W = \{\omega_0, \omega_1, \omega_2, \ldots, \omega_{n-1}\}$, $\omega_0 = \emptyset$, is a *spanning factorization* of $G$ (with word list $W$) if for every vertex $v$, the vertices $v\omega_i$ are distinct.

**Schedules**.  A schedule associated with a factorization is an assignment of a time (a label) to each occurrence of each factor in the words of $W$ such that no time is assigned more than once to a particular factor and times assigned to the factors in a single word are increasing.  The *time of a schedule* is the largest time assigned to any of the factors.  If $T$ is the total time, the schedule can be thought of as a $d \times T$ array where each row corresponds to a factor and an entry in that row indicates which occurrence of that factor has been assigned the corresponding time.  An entry in a row in the array can be empty indicating no occurrence of that factor has been assigned the given time.  The power of a spanning factorization lies in the fact that a schedule can be used to describe an algorithm for conflict free global exchange of information between the vertices of the graph.

*Theorem 3.  Suppose we have a schedule for a factorization of the graph $G$. Then the collection of directed label-increasing paths $v\omega_i$ for all $v$ and non-empty $\omega_i$ have the property that no edge in the graph is assigned the same time*

*twice. A schedule for a spanning factorization yields a time labeled directed path between every two vertices so that no edge is labeled with the same time twice.*

*Proof.* Each edge in the graph is assigned to a single one factor. Assume there is an edge in the one factor $F$ that has been assigned the same time twice. Since every occurrence of $F$ in the words in $W$ has been assigned a unique time, this can only mean that there are two different vertices $u$ and $v$ and an initial subword $\omega$ of a word in $W$ such that the edges $(u\omega, u\omega F)$ and $(v\omega, v\omega F)$ are the same edge. Then $u\omega$ and $v\omega$ must be the same vertex. Let us assume that this is the shortest $\omega$ for which this happens. The word $\omega$ cannot be empty since $u$ and $v$ are different. But then the last factor in $\omega$ must also be the same edge, a contradiction. If we start with a spanning factorization, then all the non-empty paths from $v$ are unique, there are $n-1$ of them and none of them can return to $v$ so they must reach to every other vertex in the graph.

There are some additional properties that a spanning factorization might have.

*Definition.* We say a spanning factorization is *balanced* if each factor appears nearly as often in the schedule as any other. We say the factorization is *short* if the average number of times a factor appears is the same as the theoretical lower bound $\theta$ based on the average distance between any two vertices and the number of edges. We say the factorization is optimal if it is short and balanced. A schedule $\Sigma$ is *minimum* for a spanning factorization, if it has time $\tau(\Sigma)$ equal to the theoretical minimum time for the factorization based on maximum number of times a factor appears. In mathematical terms, we can write

(i)     *balanced* if

$$\max_i |\{F_i : F_i \in W\}| = \left\lceil \sum_{i=1}^{d} |\{F_i : F_i \in W\}| / d \right\rceil ;$$

(ii)     *short* if

$$\left\lceil \sum_{i=1}^{d} |\{F_i : F_i \in W\}| / d \right\rceil = \left\lceil \sum_{k=1}^{D} \frac{k N_k}{nd} \right\rceil ,$$

here $N_k$ is the number of times the distance between two vertices is $k$ and $D$ is the diameter;

(iii)     *optimal* if

$$\max_i |\{F_i : F_i \in W\}| = \left\lceil \sum_{k=1}^{D} \frac{kN_k}{nd} \right\rceil;$$

(iv) *minimum* if

$$\tau(\Sigma) = \max_i |\{F_i : F_i \in W\}|.$$

Note that these parameters are ordered

$$\left\lceil \sum_{k=1}^{D} \frac{kN_k}{nd} \right\rceil \leq \left\lceil \sum_{i=1}^{d} |\{F_i : F_i \in W\}|/d \right\rceil \leq \max_i |\{F_i : F_i \in W\}| \leq \tau(\Sigma).$$

Creation of schedules for spanning factorizations are discussed in [6], where the following is proven.

*Fact 2.* Every diameter two spanning factorization has a minimum schedule unless the max belongs to a factor $F_i$ which is not in a word of length one and is entirely absent in words of length two in one position, either first or second. In that case, the shortest time for any schedule is one more than the theoretical minimum.

*Theorem 4.* Every Cayley graph has a short factorization.

*Proof.* This is a sketch. Take a tree $T_1$ of shortest paths from the identity of the group. The factors consist of all the edges labeled with a specific generator. The words are just the paths in $T_1$, so the factorization is short.

*Question 2.* Does every vertex symmetric graph have a factorization or even a short factorization?

**Graphs of diameter 2 with a spanning factorization.** There is one type of graph for which we can always construct an optimal transpose.

*Definition 2.* If $G$ is a directed graph and $v$ is a vertex in $G$, then for any positive integer $r$, $S_r(v)$ is the set of all vertices $u$ in $G$ for which there is a directed path from $v$ to $u$ of length less than or equal to $r$.

*Definition 3.* We call $L_r = S_r(e) - S_{r-1}(e)$ the *r layer*.

*Theorem 5. Let be a graph of diameter 2 with a spanning factorization and the elements of the 2 layer can be represented by words $x_i y_i$ in the generators. Let $M(\delta)$ be the number of times that the generator $\delta$ appears as an $x_i$ and $N(\delta)$ be the number of times that the generator $\delta$ appears as a $y_i$. If $M(\delta) + N(\delta) \leq \theta(G) - 1$ for all $\delta$, then $\tau(G) = \theta(G) = \psi(G)$.*

*Proof.* To prove this theorem, we need to show that we can find a schedule with $\theta(G)$ columns. First note that if $n_2$ is the number of vertices in the 2 layer, then

$$\sum_\delta (M(\delta) + N(\delta)) = 2n_2.$$

Since $\theta(G) = \lceil (d + 2n_2)/d \rceil = \lceil 2n_2/d \rceil + 1$, then

$$\left\lceil \sum_\delta (M(\delta) + N(\delta))/d \right\rceil = \theta(G) - 1$$

and the hypothesis says that the total number of times any generator appears is as close to the average as possible. In order to construct the schedule, we appeal to Corollary 4 in [6]. Here a generator is called a machine and a word is called a job and each generator appearing as an $x_i$ is called a type 1 step and each generator appearing as a $y_i$ is called a type 2 step.

*Corollary 4 (from [6]). Given $s_1$ jobs with one step and $s_2$ jobs with two steps, there exists an optimal schedule (in time $T = \lceil (s_1 + 2s_2)/d \rceil$) if and only if each machine is assigned at most $T$ steps and no machine is assigned exclusively steps of one type, either type 1 or type 2 from jobs with 2 steps.*

In our case, we have $s_1 = d$ and the machine (generator) $\delta$ is assigned

$$M(\delta) + N(\delta) + 1$$

steps which by hypothesis is at most $\theta(G) = T$. Since every machine is assigned one job with a single step (the members of the 1 layer), the exclusivity condition is satisfied so the hypotheses hold and the theorem guarantees the existence of the required schedule.

In addition, a ccorollary of Theorem 1 in [6] yields a more general statement.

Corollary 6. *Let be a graph of diameter 2 with a spanning factorization and the elements of the 2 layer can be represented by words $x_i y_i$ in the generators. Let $M(\delta)$ be the number of times that the generator $\delta$ appears as an $x_i$ and $N(\delta)$ be the number of times that the generator $\delta$ appears as a $y_i$. Then*

$$\tau(G) \leq 1 + \max_{\delta}(M(\delta) + N(\delta)).$$

**Comparing transpose on two networks with different parameters.** For this comparison, we need to make some model assumptions. We investigated this previously in [5] with different model assumptions. Here we create a simple model problem utilizing a transpose. The idea is that we want to compare different network architectures with the same overall cost and use the running time of a model algorithm as a way to evaluate the architecture.

First, we examine the costs. We make the simplifying assumption that memory costs are the same between architectures and only depend on the problem size. We consider only two costs, one for the processors and one for the network. We let the cost of a processor be $C_p$. For the network, we assume that cost does not depend upon how many wires are connected to a processing node but only to the overall number of wires in the network. We let the cost of one wire be $C_w$. Then in this model, the cost is

$$C = PC_p + PdC_w = P(C_p + dC_w).$$

We can write this as $\gamma = \dfrac{C}{C_w} = P\left(\dfrac{C_p}{C_w} + d\right)$.

To describe the model algorithm, we start with an $N \times N$ matrix $A$ with column vectors $v_1, v_2, \cdots, v_N$. Our algorithm also uses a given complicated function $f$ mapping $\Re^N$ to $\Re^N$. Let $F$ be the function that maps the matrix $A$ to the matrix $F(A)$ with columns $f(v_1), f(v_2), \cdots, f(v_N)$. Our algorithm is then the iterative algorithm with $M$ iterations given by

$$A_{j+1} = F(A_j^T).$$

We assume that the CPU time for this algorithm has the form $NM\alpha(N)$ for some function $\alpha$ that represents the running time of $f$. We make a simple assumption that the CPU time scales with the number of processors so that the CPU time with $P$ processors is then

$$T_p = \frac{NM\alpha(N)}{P}.$$

Given the $P$ processors and a network $G$, we distribute the columns of $A$ with $N/P$ columns in each processor. To form the transpose, each processor must send $N/P$ rows from these columns to every other processor. This means the network must send $(N/P)^2$ pieces of information to every other processor. The total communication time for the algorithm is then given by

$$T_c = M\frac{N^2}{P^2}\tau(G).$$

We can optimistically estimate $\tau$ by using ideal numbers. In the discussion, we can assume that $D$ is the average diameter of the network, that is, the sum $\sum kn_k/P$. If the network has as many processors as possible for a certain degree and diameter, it will have about $P = d^D$ processors. If we assume that the transpose is able to use all the available band width on every transpose step, then each processor has to use $D$ wires for each of the $P-1$ messages it has to transmit so the total band width requirement is $DP^2$. This yields a $\tau$ that is approximately $DP/d$ and a total communication time of

$$T_c = M\frac{N^2 D}{Pd}.$$

If we assume that $d >> C_p/C_w$, then $d^{D+1} = Pd \sim \gamma$ so the total time for the algorithm is

$$T = T_P + T_c = NM\alpha(N)\lambda^{\frac{D}{D+1}} + DN^2\lambda$$

where $\lambda = 1/\gamma < 1$.

The term $\alpha(N)$ is proportional to the number of operations in the function $f$. If $f$ were a matrix multiply, this could be anywhere from $N$ to $N^2$.

Let us examine a regime where the two terms are of comparable order in $N$. We let $\alpha(N) = \beta N$. Then we have to minimize

$$\beta\gamma^{\frac{1}{D+1}} + D.$$

Clearly then, in this regime we want $D$ to be as large as possible because that allows more processors. However, if $D$ gets too large, then $d$ will become small and the assumption that $d >> C_p/C_w$ will become invalid. In case the growth rate of $\alpha(N)$ is larger than $N$, then the $T_P$ term dominates even more strongly and again, we will want $D$ as large as possible.

We can use these results to compare two different networks with a constrained cost. First, we see that constraining cost means constraining the number of wires, that is, $Pd < \gamma = \dfrac{C}{C_w}$. Assuming both networks satisfy this constraint, the one with the largest number of processors is the better one.

**Summary**. We have discussed transpose on vertex symmetric networks. We have provided a method to compare the efficiency of transpose schemes on two different networks with a cost function based on the number processors and wires needed to complete a given algorithm in a given time.

## REFERENCES

1. Randall Dougherty and Vance Faber, "Network routing on regular directed graphs from spanning factorizations," preprint.

2. Vance Faber, "Global communication algorithms for Cayley graphs", Arxiv.org: 1305.6349.

3. Vance Faber, "Global communication algorithms for hypercubes and other Cayley coset graphs," Technical Report, LA-UR 87-3136, Los Alamos National Laboratory (1987).

4. Vance Faber, "Global sum on symmetric networks," Arxiv.org: 1201.4153.

5. Vance Faber, "Latency and diameter in sparsely populated processor interconnection networks: a time and space analysis," Technical Report, LA-UR-87-3635, Los Alamos National Laboratory (1987).

6. Vance Faber, Amanda Streib and Noah Streib, "Ideal unit-time job shop scheduling," preprint.

7. Wikipedia, en.wikipedia.org/wiki/Petersen_graph.7. Wikipedia, en.wikipedia.org/wiki/Petersen_graph.